\documentclass[12pt]{article}
\usepackage{amsxtra,amssymb,amsthm,amsmath,latexsym}

\theoremstyle{plain}

\newtheorem*{remark3.2}{Remark 3.2}
\newtheorem*{remark3.3}{Remark 3.3}
\numberwithin{equation}{section}

\def\s{{\sigma}}

\def\R{{\mathbb R}}

\def\C{{\mathbb C}}
\def\oH{\buildrel\circ\over H}
\def\oH1{\buildrel\circ\over H\kern-.02in{}^1}

\def\Im{\hbox{\,Im\,}}

\begin{document}


\title{ Inequalities for the transformation operators and applications
}

\author{
A.G. Ramm\\
LMA/CNRS, 31 Chemin Joseph Aiguier, Marseille 13402, France,\\
and Mathematics Department, Kansas State University, \\
 Manhattan, KS 66506-2602, USA\\
ramm@math.ksu.edu\\
}

\date{}

\maketitle\thispagestyle{empty}

\begin{abstract}
\footnote{Math subject classification: 34B25, 35R30, 73D25, 81F05, 81F15}
\footnote{key words: inequalities, transformation operators, inverse 
scattering, 
}
Inequalities for the transformation operator kernel
$A(x,y)$ in terms of $F$-function
are given, and vice versa. These inequalities are applied to
inverse scattering on half-line. Characterization of the scattering data 
corresponding to the usual scattering class $L_{1,1}$ of the potentials, 
to the class of compactly supported potentials, and to
the class of square integrable 
potentials is given. Invertibility of each of the steps in the inversion 
procedure is proved.

\end{abstract}


\section{Introduction}
Consider the half-line scattering problem data: 
\begin{equation} {\mathcal S} =\{ S(k), k_j, s_j, 1\leq j \leq J\},  
\end{equation}
where $ S(k)=\frac {f(-k)}{f(k)}$ is the S-matrix, $f(k)$ is the Jost
function, $f(ik_j)=0$, $\dot f(ik_j):=\frac {df(ik_j)}{dk}\neq 0$,
$k_j>0$, $s_j>0$, $J$ is a positive integer, it is equal to the number of 
negative eigenvalues of the Dirichlet operator $\ell u:=-u^{\prime 
\prime}+q(x)u$
on the half-line. The potential $q$ is real-valued throughout,
$q\in L_{1,1}:=\{q: \int_0^\infty x|q|dx<\infty \}$.
For such $q$ the scattering data ${\mathcal S}$  
have the following properties:

A) $k_j, s_j>0, \, S(-k)=\overline {S(k)}=S^{-1}(k), \, k\in \R, 
\,S(\infty)=1,$

B) $\kappa:= ind S(k):=\frac 1 {2\pi}\int_{-\infty}^\infty dlog S(k)$ is a 
nonpositive integer, 

C) $F\in L^p$, $p=1$ and $p=\infty$, $xF'\in L^1$, $L^p:=L^p(0,\infty)$.

Here 
\begin{equation} F(x):=\frac 1 {2\pi}\int_{-\infty}^\infty 
[1-S(k)]e^{ikx}dk +\sum_{j=1}^J s_j e^{-k_jx},
\end{equation}
and 
$$\kappa=-2J\quad\hbox { if } f(0)\neq 0, \quad \kappa=-2J-1 \hbox { if }
f(0)= 0.$$

The Marchenko inversion method is described in the following diagram:
\begin{equation} 
 {\mathcal S} \Rightarrow F(x) \Rightarrow A(x,y) \Rightarrow q(x),
\end{equation} 
where the step $ {\mathcal S} \Rightarrow F(x)$ is done by formula (1.2),
the step $ F(x) \Rightarrow A(x,y) $ is done by
solving the Marchenko equation:
\begin{equation}
  (I+F_x)A:=A(x,y) + \int^\infty_x A(x,t) F(t+y)\, dt = - F(x+y),
  \quad y \geq x \geq 0,
  \end{equation}
and the step $ A(x,y) \Rightarrow q(x)$ is done by the formula:
\begin{equation}
 q(x)=-2\dot A(x,x):=-2 \frac {dA(x,x)}{dx}.
\end{equation}
Our aim is to study the estimates for $A$ and $F$, which give 
a simple way to find necessary and sufficient conditions for the data
(1.1) to correspond to a $q$ from some functional class. We consider, as 
examples, 
the following classes: the usual scattering class $L_{1,1},$ for which the 
result was 
obtained earlier (\cite{M} and \cite{R} ) by a more complicated 
argument, the class of compactly 
supported potentials which are locally in  $L_{1,1}$,
and the class of square integrable potentials.
We also prove that each step in the scheme (1.3) is invertible.
In Sec.2 the estimates for $F$ and $A$ are obtained. These estimates and
their applications are the main results of the paper.
In Sec.3-6 applications to inverse scattering problem are given.
 
\section{Inequalities for $A$ and $F$}
If one wants to study the characteristic properties of the scattering data 
(1.1), that is, a necessary and sufficient condition on these data 
to guarantee that the corresponding potential belongs to a prescribed 
functional class, then conditions A) and B) are always necessary for
a real-valued $q$ to be in $L_{1,1}$, the usual class in the scattering 
theory, or other class for which the scattering theory is constructed, and 
a condition of the type C) determines actually the class of potentials 
$q$.
Conditions A) and B) are consequences of the unitarity of the 
selfadjointness of the Hamiltonian, finiteness of its negative spactrum,
and the unitarity of the $S-$matrix.
Our aim is to derive from equation (1.4) inequalities for $F$ and $A$.
This allows one do describe the set of $q$, defined by (1.5). 

Let us assume:
\begin{equation}
\sup_{y\geq x}|F(y)|:=\s_F(x)\in L^1, \quad F'\in L_{1,1}.
\end{equation}
 The function $\s_F$ is monotone decreasing, $|F(x)|\leq \s_F(x)$.
Equation (1.4) is of Fredholm type in $L^p_x:=L^p(x,\infty)$ 
$\forall x\geq 0$ and $p=1$.
The norm of the operator in (1.4) can be estimated :
\begin{equation}
||F_x||\leq \int_x^\infty \s_F(x+y)dy\leq \s_{1F}(2x), \quad 
\s_{1F}(x):=\int_x^\infty \s_F(y)dy.
\end{equation} 
Therefore (1.4) is uniquely solvable in $L^1_x$ for any $x\geq x_0$
if
\begin{equation}
\s_{1F}(2x_0)<1.
\end{equation}
This conclusion is valid for any $F$ satisfying (2.3), and conditions 
A), B), and C) are not used.
Assuming (2.3) and (2.1) and taking $x\geq x_0$, let us derive 
inequalities  for $A=A(x,y)$. Define 
$$\s_A(x):=\sup_{y\geq x}|A(x,y)|:=||A||.
$$ 
From (1.4) one gets:
$$\s_A(x)\leq \s_{F}(2x)+\s_A(x)\sup_{y\geq x}\int_x^\infty \s_F(s+y)ds 
\leq \s_{F}(2x)+\s_A(x)\s_{1F}(2x).
$$
Thus, if (2.3) holds, then
\begin{equation}
\s_A(x)\leq c\s_{F}(2x), \quad x\geq x_0.
\end{equation}
 By $c>0$  different constants depending on $x_0$ are denoted.
Let 
$$\s_{1A}(x):=||A||_1:=\int_x^\infty |A(x,s)|ds.$$
 Then (1.4) yields
$\s_{1A}(x)\leq \s_{1F}(2x)+\s_{1A}(x) \s_{1F}(2x)$. So
\begin{equation}
\s_{1A}(x)\leq c\s_{1F}(2x), \quad x\geq x_0.
\end{equation}
Differentiate (1.4) with respect to $x$ and $y$ and get:
\begin{equation}
  (I+F_x)A_x(x,y)=A(x,x)F(x+y)-F'(x+y),  
  \quad y \geq x \geq 0,
  \end{equation}
and
\begin{equation}
  A_y(x,y)+\int_x^\infty A(x,s)F'(s+y)ds=-F'(x+y),
  \quad y \geq x \geq 0.
  \end{equation}
Denote 
\begin{equation}
\s_{2F}(x):=\int_x^\infty |F'(y)|dy, \quad \s_{2F}(x)\in L^1.
\end{equation}
Then, using (2.7) and (2.4), one gets
\begin{equation}
||A_y||_1\leq  \int_x^\infty |F'(x+y)|dy+\s_{1A}(x) \sup_{s\geq 
x}\int_x^\infty |F'(s+y)|dy
\leq \s_{2F}(2x) [1+c\s_{1F}(2x)]\leq c\s_{2F}(2x),
  \end{equation}
and using (2.6) one gets:
$$||A_x||_1\leq A(x,x)\s_{1F}(2x)+\s_{2F}(2x)+||A_x||_1\s_{1F}(2x),$$ 
so
\begin{equation}
||A_x||_1\leq c[\s_{2F}(2x)+\s_{1F}(2x)\s_{F}(2x)].
  \end{equation}

Let $y=x$ in (1.4), then differentiate (1.4) with respect to $x$
and get:
\begin{equation}
\dot A(x,x)= -2F'(2x)+A(x,x)F(2x)- \int_x^\infty 
A_x(x,s)F(x+s)ds -\int_x^\infty A(x,s)F'(s+x)ds.
  \end{equation}
From (2.4), (2.5), (2.10) and (2.11) one gets:
\begin{equation}  
|\dot A(x,x)|\leq 
2|F'(2x)|+c\s^2_F(2x)+c\s_F(2x)[\s_{2F}(2x)+\s_{1F}(2x)\s_{F}(2x)]
+c\s_F(2x)\s_{2F}(2x).
  \end{equation}
Thus,
\begin{equation}
x|\dot A(x,x)|\in L^1, 
  \end{equation}
provided that $xF'(2x)\in L^1, \, x\s^2_F(2x)\in L^1,$ and
 $ x\s_F(2x)\s_{2F}(2x)\in L^1.$
Assumption (2.1) implies $xF'(2x)\in L^1$. If $\s_{F}(2x)\in L^1$,
and $\s_{F}(2x)>0$ decreases monotonically, then $x\s_{F}(x)\to 0$
as $x\to \infty$. Thus $x\s^2_F(2x)\in L^1,$ and
$ \s_{2F}(2x)\in L^1$  because $\int_0^\infty dx\int_x^\infty |F'(y)|dy=
\int_0^\infty |F'(y)|ydy<\infty$, due to (2.1).
Thus, (2.1) implies (2.4), (2.5), (2.8), (2.9), and (2.12), while (2.12) 
and (1.5) imply $q\in \tilde L_{1,1}$ where $\tilde 
L_{1,1}=\{q: q={\overline q},\, 
\int^\infty _{x_0}x|q(x)|dx<\infty\}$, and $x_0\geq 0$ satisfies (2.3).

Let us assume now that (2.4), (2.5), (2.9), and (2.10) hold, where 
$\s _F\in L^1$ and $\s_{2F}\in L^1$ are some positive monotone decaying 
functions (which have nothing to do now with the function $F$,
solving equation (1.4), and derive 
estimates for this function $F$. let us rewrite (1.4) as:
\begin{equation}
F(x+y)+\int^\infty_x A(x,s)F(s+y)ds=-A(x,y), \qquad y\geq x\geq 0.
\end{equation}
Let $x+y=z, s+y=v$. Then,
\begin{equation}
F(z)+\int^\infty_z A(x,v+x-z)F(v)dv=-A(x,z-x),\qquad z\geq 2x.
\end{equation}
From (2.15) one gets:
$$\s_F(2x)\leq\s_A(x)+
\s_F(2x)\sup_{z\geq 2x}\int^\infty_z |A(x,v+x-z)|dv\leq\s_A(x)+\s_F(2x)\, 
||A||_1.
$$

Thus, using (2.5) and (2.3), one obtains:
\begin{equation}
\s_F (2x)\leq c\s_A(x).
\end{equation}
Also from (2.15) it follows that:
\begin{equation}
\begin{array}{ll}
\s_{1F}(2x)&:=||F||_1:=\int^\infty_{2x}|F(v)|dv \\
&\leq\int^\infty_{2x}
|A(x,z-x)|dz+\int^\infty_{2x}\int^\infty_z |A(x,v+x-z)| |F(v)|dvdz \\
&\leq
||A||_1 +||F||_1 ||A||_1, \\
&\text {so} \\
& \s_{1F}(2x)\leq c\s_{1A}(x).
\end{array}
\end{equation}
From (2.6) one gets:
\begin{equation}
\int^\infty_x |F'(x+y)|dy=\s_{2F}(2x)\leq c\s_A (x) \s_{1A}(x)+||A_x ||+c 
||A_x||_1 \s_{1A}(x).
\end{equation}
Let us summarize the results:

{\bf Theorem 2.1.} {\it If $x\geq x_0$ and (2.1) hold, then one has:
\begin {equation}
\begin{array}{ll}
\s_A(x)\leq c\s_F (2x),\quad \s_{1A}(x)\leq c\s_{1F}(2x),\quad 
||A_y||_1\leq\s_{2F}(2x)(1+c\s_{1F}(2x)),\\
||A_x||_1\leq 
c[\s_{2F}(2x)+\s_{1F}(2x)\s_F (2x)].
\end{array}
\end{equation}
Conversely, if $x\geq x_0$ and
\begin{equation}
\s_{A}(x)+\s_{1A}(x) +||A_x ||_1 +||A_y||_1 <\infty,
\end{equation}
then
\begin{equation}
\begin{array}{ll}
\s_F (2x)\leq c\s_{A}(x),\quad \s_{1F}(2x)\leq c\s_{1A}(x),&\\ 
\s_{2F}(x)\leq 
c[\s_{A}(x)\s_{1A}(x) +||A_x||_1 (1+\s_{1A}(x))].
\end{array}
\end{equation}
}
In Sec. 3 we replace the assumption $x\geq x_0>0$ by $x\geq 0$. The 
argument 
in this case is based on the Fredholm alternative.

\section{ Applications}

First, let us give necessary and sufficient conditions on 
${\mathcal S}$ for $q$ to 
be in $L_{1,1}$. These conditions are known \cite{M}, \cite{R}
and \cite{R1}, but we give a short new argument using some ideas
from \cite{R1}. We 
assume throughout that conditions A), B), 
and C) hold. These conditions are known to be necessary for $q\in 
L_{1,1}$. Indeed, conditions A) and B) are obvious, and C) is proved in 
Theorems 
2.1 and 3.3. Conditions A), B), and C) are also sufficient for $q\in 
L_{1,1}$. Indeed if they 
hold, then we prove that equation (1.4) has a unique solution in $L^1_x$ 
for all $x\geq 0$. This is a known fact \cite{M}, but we give a (new) 
proof because it is short. This proof combines some ideas from \cite{M} 
and \cite{R1}.

{\bf Theorem 3.1.}  {\it If A), B), and C) hold, then (1.4) has a 
solution in 
$L^1_x$ for any $x\geq 0$ and this solution is unique.}

{\bf Proof.} Since $F_x$ is compact in $L^1_x,\, \forall x\geq 0$, by 
the Fredholm 
alternative it is sufficient to prove that 
\begin{equation}
(I+F_x)h=0,\quad h\in L^1_x, 
\end{equation}
implies $h=0$. Let us prove it for $x=0$. The proof is similar for 
$x>0$. If $h\in L^1$, then $h\in L^\infty$ because $||h||_\infty \leq 
||h||_{L^1}\s_F(0)$. If $h\in L^1\cap L^\infty$, then $h\in 
L^2$ because$ ||h||_{L^2}^2\leq ||h||_{L^\infty}||h||_{L^1}$. Thus,
if $h\in L^1$ and solves (3.1), then $h\in 
L^2\cap L^1\cap L^\infty$.

Denote $\tilde h=\int^\infty_0 h(x)e^{ikx}dx,\, h\in L^2$. Then, 
\begin{equation}
\int^\infty_{-\infty}\tilde h^2 dk=0.
\end{equation}
Since $F(x)$ is real-valued, one can assume $h$  real-valued. 
One has, using Parseval's equation:
$$
0=((I+F_0)h,h)=\frac 1{2\pi} ||\tilde h||^2+\frac 
1{2\pi}\int^\infty_{-\infty}[1-S(k)]\tilde 
h^2(k)dk+\sum^J_{j=1}s_jh^2_j,\quad h_j:=\int^\infty_0e^{-k_jx}h(x)dx.
$$
Thus, using (3.2), one gets 
$$h_j=0,\, 1\leq j\leq J,\quad (\tilde h,\tilde h)= 
(S(k)\tilde h,\, \tilde h(-k)),
$$ where we have used 
real-valuedness of $h$, i.e. $ \tilde h(-k)=\tilde h(k),\forall 
k\in R$.

Thus, $(\tilde h, \tilde h)=(\tilde h,S(-k)\tilde h(-k))$, where A) was 
used. Since $||S(-k)||=1$, one has $||\tilde h||^2=\vert (\tilde 
h,S(-k)\tilde h(-k))\vert\leq ||\tilde h||^2$, so the equality sign is 
attained in the Cauchy inequality. Therefore, $\tilde h(k)=S(-k)\tilde 
h(-k)$.

By condition B), the theory of Riemann problem (see [1]) guarantees 
existence and uniqueness of an analytic in $\C_+:=\{k: \Im k>0\}$ function 
$f(k):=f_+(k),\, f(ik_j)=0,\, \dot f(ik_j)\neq 0,\, 1\leq j\leq J, 
\, f(\infty)=1$, such that 
\begin{equation}
f_+(k)=S(-k)f_-(k), \quad k\in\R,
\end{equation}
and $f_-(k)=f(-k)$ is analytic in $\C_- :=\{k:Imk<0\},\, f_-(\infty)=1$ 
in $\C_-,\, f_-(-ik_j)=0,\, \dot f_-(-ik_j)\neq 0$. Here the property 
$S(-k)=S^{-1}(k),\, \forall k\in \R$ is used.

One has
$$
\psi (k):=\frac{\tilde h(k)}{f(k)}=\frac{\tilde h(-k)}{f(-k},\quad k\in 
\R,\quad h_j=\tilde h(ik_j)=0,\quad 1\leq j\leq J.
$$
The function $\psi (k)$ is analytic in $\C _+$ and $\psi (-k)$ is analytic 
in $\C _-$, they agree on $\R$, so $\psi(k)$ is analytic in $\C$. Since 
$f(\infty )=1$ and $\tilde h(\infty )=0$, it follows that $\psi\equiv 0$.

Thus, $\tilde h=0$ and, consequently, $h(x)=0$, as claimed.
Theorem 3.1 is proved. $\Box$

The unique solution to equation (1.4) satisfies the estimates given in 
Theorem 2.1. In the proof of Theorem 2.1 the estimate $x|\dot A(x,x)|\in 
L^1(x_0,\infty )$ was established. So, by (1.5), $xq\in L^1(x_0,\infty)$.

The method developed in Sec.2 gives accurate information about the 
behavior of $q$ near infinity. An immediate consequence of Theorems 2.1 
and 3.1 is:

{\bf Theorem 3.2.} {\it If A), B), and C) hold, then $q,$ obtained by the 
scheme (1.3), belongs to $L_{1,1}(x_0,\infty )$.}

Investigation of the behavior of $q(x)$ on $(0,x_0)$ requires additional 
argument. Instead of using the contraction mapping principle and 
inequalities, as in Sec. 2, one has to use the Fredholm theorem, which 
says that $||(I+F_x)^{-1}||\leq c$ for any $x\geq 0$ where the operator 
norm is for $F_x$ acting in $L^p_x$, $ p=1$ and $p=\infty$, and the 
constant $c$ does not depend on $x\geq 0$.

Such an analysis yields:

{\bf Theorem 3.3.} {\it If and only if A), B), and C) hold, then $q\in 
L_{1,1}$.}

{\bf Proof.} It is sufficient to check that Theorem 2.1 holds with $x\geq 
0$ replacing $x\geq x_0$. To get (2.4) with $x_0=0$, one uses (1.4) and 
the 
estimate:
\begin{equation}
||A(x,y)||\leq ||(I+F_x)^{-1}||||F(x+y)||\leq c\s_F(2x), \quad ||\cdot||=
\sup_{y\geq x}|\cdot|,\, x\geq 0,
\end{equation}
where the constant $c>0$ does not depend on $x$. Similarly:
\begin{equation}
||A(x,y)||_1\leq c\sup_{s\geq x} \int_x^\infty |F(s+y)|dy\leq 
c\s_{1F} (2x), \quad x\geq 0.
\end{equation}
From (2.6) one gets:
\begin{equation}
\begin{array} {ll}
||A_x(x,y)||_1\leq c[||F'(x+y)||_1 +A(x,x)||F(x+y)||_1]&\\\leq 
c\s_{2F}(2x)
+c\s_{F}(2x)\s_{1F}(2x),
\quad x\geq 0. 
\end{array}
\end{equation}
From (2.7) one gets:
\begin{equation}
||A_y(x,y)||_1\leq c[\s_{2F}(2x) 
+\s_{1F}(2x)\s_{2F}(2x)]\leq \s_{2F}(2x).
\end{equation}
Similarly, from (2.11) and (3.3)-(3.6) one gets 
(2.12). Then one checks (2.13) as in the proof of Theorem 2.1.
Consequently Theorem 2.1 holds with $x_0=0$.
Theorem 3.3 is proved. $\Box$.

\section{ Compactly supported potentials}

In this Section necessary and sufficient conditions are
given for $q\in L^a_{1,1}:=\{ q: q={\overline q}, q=0  \hbox { if }
x>a, \int_0^a x|q|dx<\infty\}$. Recall that the Jost solution is:
\begin{equation}
f(x,k)=e^{ikx}+\int_x^\infty A(x,y)e^{iky}dy,\quad f(0,k):=f(k).
\end{equation}

{\bf Lemma 4.1.} {\it If $q\in L^a_{1,1}$, then $f(x,k)=e^{ikx}$ for 
$x>a$,
$A(x,y)=0$ for $y\geq x\geq a$,  $F(x+y)=0$ for $y\geq x\geq a$
(cf (1.4)), and $F(x)=0$ for $x\geq 2a$.}

Thus, (1.4) with $x=0$ yields $A(0,y):=A(y)=0$ for $x\geq 2a$.
The Jost function 
\begin{equation}
f(k)=1+\int_0^{2a} A(y)e^{iky}dy,\quad A(y)\in W^{1,1}(0,a),
\end{equation}
is an entire function of exponential type $\leq 2a$, that is,
$|f(k)|\leq ce^{2a|k|}$, $k\in \C$, and $S(k)=f(-k)/f(k)$ is a 
meromorphic function in $\C$. In (4.2) $W^{l,p}$ is the Sobolev space,
and the inclusion (4.2) follows from Theorem 2.1.

Let us formulate the assumption D):

{\it D) the Jost function $f(k)$ is an entire function of exponential type 
$\leq 2a$.}

{\bf Theorem 4.1.} {\it Assume A),B), C) and D). Then $q\in L^a_{1,1}$. 
Conversely, if  $q\in L^a_{1,1}$, then A),B), C) and D) hold.}

{\bf Proof.} {\it Necessity.} If $q\in L_{1,1}$, then A), B) and C) hold 
by
Theorem 3.3, and D) is proved in Lemma 4.1. The necessity is proved. 

{\it Sufficiency.} If  A), B) and C) hold, then $q\in L_{1,1}$. One has 
to
prove that $q=0$ for $x>a$. If D) holds, then from the proof of Lemma 4.1 
it follows that $A(y)=0$  for $y\geq 2a$. 

{\it We claim that $F(x)=0$ for $x\geq 2a$}. 

If this is proved, then (1.4) yields
$A(x,y)=0$ for $y\geq x\geq a$, and so $q=0$ for $x>a$ by (1.5).

Let us prove the claim. 

Take $x>2a$ in (1.2). The function $1-S(k)$
is analytic in $\C_+$ except for $J$ simple poles at the points
$ik_j$. If $x>2a$ then one can use the Jordan lemma and residue 
theorem and get:
\begin{equation}
F_S(x)=\frac 1 {2\pi} \int_{-\infty}^\infty 
[1-S(k)]e^{ikx}dk=-i\sum^J_{j=1}\frac{f(-ik_j)}{\dot 
f(ik_j)}e^{-k_jx},\quad 
x>2a.
\end{equation}
Since $f(k)$ is entire, the Wronskian formula 
$$f^\prime (0,k)f(-k)-f^\prime (0,-k)f(k)=2ik
$$ 
is valid on $\C$, and at $k=ik_j$ it yields: 
$$f^\prime(0,ik_j)f(-ik_j)=-2k_j,
$$
because $f(ik_j)=0$.
 This and (4.3) yield
$$
F_s(x)=\sum^J_{j=1}\frac{2ik_j}{f^\prime(0,ik_j)\dot 
f(ik_j)}e^{-k_jx}=-\sum ^J_{j=1}s_je^{-k_jx}=-F_d(x),\quad x>2a.
$$
Thus, $F(x)=F_s(x)+F_d(x)=0$ for $x>2a$. The sufficiency is proved.

Theorem 4.1 is proved.$\Box$

In \cite{M} a condition on ${\mathcal S}$, which guarantees that $q=0$ for 
$x>a$, is given 
under the assumption that there is no discrete spectrum, that is $F=F_s$.

\section{ Square integrable potentials}

Let us introduce conditions (5.1) - (5.3):
\begin{equation}
2ik[f(k)-1+\frac Q{2ik}]\in L^2(\R):=L^2,\quad Q:=\int^\infty _0qds,
\end{equation}
\begin{equation}
k[1-S(k)+\frac Q{ik}]\in L^2,
\end{equation}
\begin{equation}
k[|f(k)|^2-1]\in L^2.
\end{equation}
{\bf Theorem 5.1.} {\it If A), B), C), and any one of the conditions
(5.1)-(5.3) hold, then $q\in L^2(\R)$.}

{\bf Proof.} We refer to \cite{R} for the proof. $\Box$

\section{ Invertibility of the steps in the inversion procedure} 
We assume A), B), and C) and prove:

{\bf Theorem 6.1.} {The steps in (1.3) are invertible:}
\begin{equation}
{\mathcal S}\Longleftrightarrow F\Longleftrightarrow A\Longleftrightarrow 
q.
\end{equation}

{\bf Proof.} 
1. Step ${\mathcal S}\Rightarrow F$ is done by formula (1.2). Step 
$F\Rightarrow {\mathcal S}$ is done by taking $x\rightarrow -\infty$ in (1.2). 
The asymptotics of $F(x)$, as $x\rightarrow -\infty$, yields 
$J, s_j, k_j,\quad 1\leq j\leq J$, that is, $F_d(x)$. Then $F_s=F-F_d$ is 
calculated, and $1-S(k)$ is calculated by taking the inverse Fourier 
transform 
of $F_s(x)$. Thus, $\mathcal S$ is found.

2. Step $F\Rightarrow A$ is done by solving (1.4), which has one and only 
one solution in $L^1_x$ for any $x\geq 0$ by Theorem 3.1. Step 
$A\Rightarrow F$ is done by solving equation (1.4) for $F$. Let $x+y=z$ 
and $s+y=v$. Write (1.4) as
\begin{equation}
(I+B)F:=F(z)+\int_z^\infty A(x,v+x-z)F(v)dv=-A(x,z-x),\quad z\geq 
2x\geq 0.
\end{equation}
The norm of the integral operator $B$ in $L^1_{2x}$ is estimated as 
follows:
\begin{equation}
\begin{array}{ll}
||B||\leq\sup_{v>0}\int^v_0|A(x,v+x-z)|dz\leq 
c\sup_{v>0}\int^v_0\sigma 
(x+\frac{v-z}2)dz&\leq\\
2\int^\infty_0\sigma(x+
w)dw=2\int^\infty_x\sigma (t)dt,
\end{array}
\end{equation}
where the known estimate \cite{M} was used:
 $|A(x,y)|\leq c\sigma (\frac{x+y}2),\quad \sigma (x):=\int ^\infty 
_x|q|dt$. It follows from (6.3) that $||B||<1$ if $x>x_0$, where $x_0$ is 
large enough. Indeed, $\int^\infty_x\sigma(s)ds\rightarrow 0$ as 
$x\rightarrow\infty$ if $q\in L_{1,1}$. Therefore, for $x>x_0$ equation 
(6.2) is uniquely solvable in $L^1_{2x_0}$ by the contraction mapping 
principle.

3. Step $A\Rightarrow q$ is done by formula (1.5). Step $q\Rightarrow A$ 
is 
done by solving the known Volterra equation (see \cite{M} or \cite{R}):
\begin {equation}
A(x,y)=\frac{1}{2} \int^\infty_{\frac{x+y}{2}} q(t) \,dt
  + \int^\infty_{\frac{x+y}{2}}
  ds\int^{\frac{y-x}{2}}_0\, dt q(s-t)A(s-t,s+t).
\end{equation}
Theorem 6.1 is proved.$\Box$

Note that Theorem 6.1 implies that if one starts with a $q\in L_{1,1}$, 
computes the scattering data (1.1) corresponding to this $q$, and uses the 
inversion scheme (1.3), then the potential obtained by the formula (1.5) 
is equal to the original potential $q$.

If $F(z)$ is known for $x\geq 2x_0$, then (6.2) can be written as a 
Volterra equation with a finite region of integration.
\begin{equation}
F(z)+\int^{2x_0}_zA(x,v+x-z)F(v)dv=-A(x,z-x)-\int^\infty_{2x_0}
A(x,v+x-z)F(v)dv,
\end{equation}
where the right-hand side in (6.5) is known.
This Volterra integral equation on the interval $z\in (0,2x_0)$ is 
uniquely solvable by iterations. Thus, $F(z)$ is uniquely determined on 
$(0,2x_0)$, and, consequently, on $(0,\infty)$.

\end{document}